\documentclass[12pt]{article}
\input epsf
\usepackage{amsmath}
\usepackage{amssymb}
\usepackage{amscd}
\usepackage{epsfig}

\def\no{\noindent}
\newenvironment{conjecture}{\noindent{\bf Conjecture}\em}{}
\newtheorem{lemma}[equation]{Lemma}
\newtheorem{lem}[equation]{Lemma}
\newtheorem{corollary}[equation]{Corollary}
\newtheorem{proposition}[equation]{Proposition}
\newtheorem{theorem}[equation]{Theorem}
\newtheorem{definition}{Definition}

\newcommand{\lra}{\longrightarrow}
\newcommand{\ra}{\rightarrow}
\newcommand{\restr}{\mbox{\Large \(|\)\normalsize}}

\newcommand{\R}{\mathbb R}

\newcommand{\Z}{\mathbb Z}

\renewcommand{\H}{\mathbb H}

\newcommand{\acts}{\curvearrowright}

\newcommand{\qed}{{\unskip\nobreak\hfil
        \penalty50\hskip1em\hbox{}\nobreak\hfil
        $\square$\parfillskip=0pt\finalhyphendemerits=0 \par}}
\newcommand{\proof}{\no{\em Proof.\ }}

\def\D{\partial}

\def\de{\delta}

\def\eps{\epsilon}
\def\ga{\gamma}
\def\be{\beta}
\def\Ga{\Gamma}

\def\lang{\langle}
\def\<{\lang}
\def\>{\rangle}
\def\lra{\longrightarrow}

\def\ol{\overline}

\def\ra{\rightarrow}

\def\si{\sigma}
\def\Si{\Sigma}

\def\geo{\partial_{\infty}}

\def\defeq{:=}
\def\hook{\hookrightarrow}

\renewcommand{\S}{{\mathcal S}}
\newcommand{\trip}{\partial^3}
\newcommand{\dub}{\partial^2_{\infty}}
\newcommand{\db}{\bar\partial^2_{\infty}}
\newcommand{\dxt}{\widetilde{DX}}
\newcommand{\geodx}{\geo\widetilde{DX}}
\newcommand{\normal}{\vartriangleleft}

\begin{document}

\title{Hyperbolic groups with 1-dimensional boundary}
\author{Michael Kapovich\thanks{Supported by NSF grant
DMS-96-26633}\\
Bruce Kleiner\thanks{Supported by a Sloan Foundation
Fellowship, and NSF grants
DMS-95-05175, DMS-96-26911, DMS-9022140.}}
\date{June 10, 1998}
\maketitle

\begin{abstract}
If a torsion-free hyperbolic group $G$ has 1-dimensional
boundary $\geo G$, then $\geo G$ is a Menger curve
or a Sierpinski carpet provided $G$ does not split over
a cyclic group.  When $\geo G$ is a Sierpinski carpet we show that $G$
is a quasi-convex subgroup of a $3$-dimensional hyperbolic  
Poincar\'e duality group. We also 
construct a ``topologically rigid'' hyperbolic group $G$: any 
homeomorphism of $\geo G$ is induced by an element of $G$. 
\end{abstract}

\setcounter{section}{1}
\setcounter{subsection}{0}

\subsection{Introduction}

\no We recall that the boundary $\geo X$ of a locally compact 
Gromov hyperbolic space $X$ is a compact metrizable topological space. 
Brian Bowditch observed that any compact metrizable space $Z$ arises this
way: view the unit ball $B$ in Hilbert space as the Poincar\'e model of 
infinite dimensional hyperbolic space, topologically embed $Z$
in the boundary of $B$, and then take the convex hull $CH(Z)$ to 
get a locally compact Gromov hyperbolic space with $\geo CH(Z)=Z$. 
 On the other hand when 
$X$ is the Cayley graph of a Gromov hyperbolic group $G$  
then the topology of $\geo X\simeq\geo G$ is quite restricted. It is known that $\geo G$ 
is finite dimensional, and either perfect,  empty, or a
two element set (in the last two cases the group $G$ is {\bf elementary}). 
It was shown recently by Bowditch and Swarup \cite{Bowditch2, Swarup} 
that if $\geo G$ is connected then it does not have global cut-points, 
and thus is locally connected according to \cite{BM}. The boundary of $G$ 
necessarily has a ``large'' group of homeomorphisms: if $G$ is nonelementary 
then its action on $\geo G$ is minimal, and 
$G$ acts on $\geo G$ as a discrete uniform convergence group. It 
turns out that the last property gives a dynamical characterization 
of boundaries of hyperbolic groups, according to a theorem of Bowditch 
\cite{Bowditch3}: if $Z$ is a compact metrizable space with $|Z|\geq 3$ and $G
\subset Homeo(Z)$ is a discrete uniform convergence subgroup, then $G$ is 
hyperbolic and $Z$ is $G$-equivariantly homeomorphic to $\geo G$. 

\medskip
There are two questions which arise naturally:

\medskip
\no {\bf Question A.}  {\em Which topological spaces are boundaries of 
hyperbolic groups?}

\medskip
\no {\bf Question B.} {\em Given a topological space $Z$, which 
hyperbolic groups have $Z$ as the boundary?} 

\medskip
Regarding question A, all finite-dimensional topological spheres and some 
homology spheres \cite{Davis}, the Sierpinski carpet and the 
Menger   curve \cite{Benakli} arise as boundaries of hyperbolic groups. 
Moreover, according to Gromov and Champetier \cite{Champetier}, 
``generic'' finitely presentable groups are hyperbolic and have the Menger 
curve as the boundary. On the other hand, as was noticed by Bestvina, it 
is unknown if higher-dimensional Universal
Menger compacta \cite{Bestvina} appear as boundaries of 
hyperbolic groups (Dranishnikov can construct hyperbolic groups 
with boundary homeomorphic to the 2-dimensional Menger compactum, \cite{Dr}).  

Considerably less is known about the Question B. If $\geo G$ is 
zero-dimensional then $G$ is a virtually free group
  \cite{Stallings,Gromov,Swiss}. 
 Recently it was proven 
in \cite{Gabai, Casson, Tukia} that any hyperbolic group whose boundary 
is homeomorphic to $S^1$  acts discretely, cocompactly, and isometrically
on the hyperbolic 
plane. We call such group $G$ {\bf virtually Fuchsian}.
The case when $\geo G \simeq S^2$ is a difficult open problem:

\begin{conjecture} (J. Cannon).  If $G$ is a hyperbolic group whose 
boundary is homeomorphic to $S^2$, then $G$ acts isometrically and properly 
discontinuously on hyperbolic 3-space $\H^3$.
\end{conjecture}

One can also answer Question B  for 
 topologically rigid hyperbolic groups:

\begin{definition} A hyperbolic group $G$ is said to be {\bf topologically rigid}
if every homeomorphism $f: \geo G \to \geo G$ is induced by an element of $G$.
\end{definition} 

\medskip 
\no
If $G'$ is a  hyperbolic group whose boundary is homeomorphic
to the boundary of a topologically rigid hyperbolic group $G$, then there
is a finite normal subgroup $N\normal G'$ so that $G'/N$ embeds in 
$G$ as a finite index subgroup. We construct  a topologically rigid
group in the section \ref{toprigid}. 

\medskip
Most of the results of our paper concern hyperbolic groups with 
one-dimensional boundary. 

\begin{theorem}
\label{s1mengersierpinski}
Let $G$ be a hyperbolic group which does not split over
a finite or virtually cyclic subgroup, and suppose $\geo G$
is 1-dimensional.  Then one of the following holds (see section
\ref{prelim} for definitions):

1. $\geo G$ is a Menger curve.

2. $\geo G$ is a Sierpinski carpet.

3. $\geo G$ is homeomorphic to $S^1$ and 
$G$ maps onto a Schwartz triangle group with finite kernel.
\end{theorem}

\medskip
It is probably impossible to classify hyperbolic 
groups whose boundaries are homeomorphic to the Menger curve
(since this is the ``generic'' case), however 
it appears that a meaningful study is possible in the case of  hyperbolic 
groups whose boundaries are homeomorphic to the Sierpinski carpet. 
Recall that the Sierpinski carpet ${\cal S}$ has canonical collection of 
{\bf peripheral circles} (see section \ref{prelim}).    

\begin{theorem}
\label{Doubletheorem}
Suppose that $\geo G \cong {\cal S}$. Then:

1. There are only finitely many $G$-orbits of peripheral circles. 

2. The stabilizer of each peripheral circle $C$ is a quasi-convex virtually Fuchsian 
group which acts on $C$ as a uniform convergence group. We call these 
subgroups {\bf peripheral subgroups of} $G$. 

3. If we ``double'' $G$ along the collection of peripheral subgroups 
using  amalgamated free product and iterated HNN-extension
(see section \ref{Doubling}), then the result 
is a hyperbolic group $\hat{G}$ which contains $G$ as a quasiconvex subgroup. 

4. The boundary of $\hat{G}$ is homeomorphic to $S^2$. Hence
by \cite{BM}, \cite{Best}, $\hat{G}$ is a $3$-dimensional
Poincar\'e duality group in the torsion-free case. 

5. When $G$ is torsion free, then $(G;H_1,\ldots,H_k)$ is a 3-dimensional 
relative Poincar\'e duality pair (see \cite{DD} for the definition).
\end{theorem}

Known examples are consistent with the following:

\begin{conjecture}
Let $G$ be a hyperbolic group Sierpinski carpet boundary.
Then $G$ acts discretely, cocompactly, and isometrically
on a convex subset of $\H^3$ with nonempty totally geodesic boundary.
\end{conjecture}

There is now some evidence for this conjecture.   It would
follow from a positive
solution of Cannon's conjecture together with Theorem
\ref{Doubletheorem}, see section \ref{Doubling}.  Alternately,
in the torsion-free case, if one could
show that (hyperbolic) $3$-dimensional Poincar\'e duality  groups are $3$-manifold
groups, then Thurston's Haken uniformization theorem 
could be applied to an irreducible $3$-manifold with 
fundamental   group isomorphic to the group $\hat{G}$ produced
in Theorem \ref{Doubletheorem}.  Under extra conditions (such as coherence
and the existence of a nontrivial splitting)
one can show that a $3$-dimensional Poincar\'e duality group
is a $3$-manifold group, \cite{KK}.

The conjecture above leads one to ask 
which hyperbolic groups have planar boundary.  Concretely,
one may ask if a torsion-free hyperbolic group with planar
boundary has a finite index subgroup subgroup isomorphic
to a discrete convex cocompact subgroup of $Isom(\H^3)$.
Here is a cautionary example which shows that in general it is necessary
to pass to a finite index subgroup: if one takes a surface of genus $1$
with two boundary components and glues one boundary circle to 
the other by a degree $2$ map, then the fundamental group $G$
of the resulting complex enjoys the following properties:

1.  $G$ is torsion-free and  hyperbolic.

2.  $G$ contains a finite index subgroup which isomorphic to 
a discrete, convex cocompact subgroup of $Isom(\H^3)$. In particular, 
the boundary of $G$ is planar. 

3. $G$ is not a $3$-manifold group.

\subsection{Preliminaries}
\label{prelim}

\no
{\bf Properties of hyperbolic groups and spaces.}
For a proof of the following properties of hyperbolic groups,
we refer the reader to \cite{Gromov,Many,Swiss,Bowditch3}.

Let $G$ be a nonelementary Gromov hyperbolic group, and suppose
$G$ acts discretely and cocompactly on a locally compact geodesic
metric space $X$.
Then the boundary of $X$ is a compact metrizable space $\geo X$
on which $Isom(X)$ acts by homeomorphisms.  For any $f\in Isom(X)$,
we denote the corresponding homeomorphism of $\geo X$ by
$\geo f$.  The  action of $G$ on $\geo X$ is minimal,
i.e. the $G$-orbit of every point is dense in $\geo X$.  Let
$\dub X\defeq \geo X\times \geo X - Diag$ be the space of distinct
pairs in $\geo X$.  Then the set of pairs of points $(x,y)\in\dub X$ 
which are fixed by an infinite cyclic subgroup of 
$G$ is dense in $\dub X$. We let $\db X:= \dub X/(x,y)\sim (y,x)$. 

 The group $G$ acts  
cocompactly and properly discontinuously on $\trip X
\defeq\{ (x,y,z)\in (\geo X)^3\mid \mbox{$x,\,y,\,z$ distinct}\}$.
There is a natural topology on $X\cup\geo X$  which is a $G$-invariant
compactification of $X$, and this is compatible with the 
topology on $\geo X$.

Recall that a subset $S$ of a geodesic metric space
is $C$-quasi-convex if every geodesic segment with
endpoints in $S$ is contained in the $C$-tubular
neighborhood of $S$. Quasi-convex subsets of
$\de$-hyperbolic metric spaces satisfy
a {\bf visibility property} (cf. \cite{EberONeill}):

\begin{quotation}
{\em Given $R,\,C,\,\de\in (0,\infty)$ there is an $R'$
with the following property (we may take $R'=R+10\de$).
  If $X$ is a $\de$-hyperbolic
metric space, $Y\subset X$ is $C$-quasi-convex, and 
$x\in X$ satisfies $d(x,Y)\geq R'$, then given any
two unit speed geodesics $\ga_1,\,\ga_2$ starting at 
$x$ and ending in $Y$, and any $t\in [0,R]$
we have $d(\ga_1(t),Im(\ga_2))<\de$ and $d(\ga_2(t),Im(\ga_1))<\de$.}
\end{quotation}

As a consequence of the visibility property, if $Y_k\subset X$
is a sequence of $C$-quasi-convex subsets of a $\de$-hyperbolic
space $X$, and $d(x,Y_k)\ra\infty$ as $k\ra \infty$, then
a subsequence of $Y_k$'s converges to a single point $\xi\in\geo X$.

\medskip
\no
{\bf Sierpinski carpets and Menger curves.}
The classical construction of a Sierpinski carpet is 
analogous to the construction of a Cantor set: 
start with the unit square in the plane, subdivide it into 
nine equal subsquares, remove the middle open square, and
then repeat this procedure inductively on the remaining squares.
If we take a sequence $D_i\subset S^2$ of disjoint closed
2-disks whose union is dense in $S^2$ so that $Diam(D_i)\ra 0$ as 
$i\ra \infty$, then $S^2-\cup_iInterior(D_i)$ is a Sierpinski carpet; 
moreover any Sierpinski carpet embedded in $S^2$ is obtained in this
way \cite{Whyburn}.  Sierpinski carpets
can also be characterized as follows \cite{Whyburn}:
a compact, 1-dimensional, planar, connected, locally connected
space with no local cut points is a Sierpinski carpet.

We will use a few topological
properties of Sierpinski carpets  ${\cal S}$:

1.  There is a unique embedding of $\S$ in $S^2$
up to post-composition with a  homeomorphism of
$S^2$.

2.  There is a countable collection ${\cal C}$ of ``peripheral
circles'' in $\S$, which are precisely the nonseparating
topological circles in $\S$.

3.  Given any metric $d$ on $\S$ and any number $D>0$,
there are only 
finitely many peripheral circles in $\S$ of
diameter $>D$.

The Menger curve may be constructed as follows. Start with the 
unit cube $I^3$ in $\R^3$.  Consider the orthogonal projections
$\pi_{ij}:I^3\ra F_{ij}$ of the unit cube onto the $ij$
coordinate square, and let ${\cal S}_{ij}\subset F_{ij}$
be the Sierpinski carpet as constructed above.
The Menger curve is the intersection 
$\cap_{i<j}\pi_{ij}^{-1}({\cal S}_{ij})$.  
   The Menger curve is universal among
all compact metrizable 1-dimensional spaces: any such
space can topologically embedded in the Menger curve.
By \cite{Anderson1,Anderson2}, a compact, metrizable, connected,
locally connected, 1-dimensional space is a Menger curve
provided it has no local cut points, and no nonempty
open subset is planar.

\subsection{Proof of Theorem \ref{s1mengersierpinski}}

The fact that $G$ does not split over a finite group implies
\cite{Stallings} that $G$ is one-ended, and $\geo G$ is 
connected.
Recall that by the results of \cite{BM,Bowditch2,Swarup}, the
boundary of a  one-ended
hyperbolic group is locally connected and has no global cut points;  
furthermore,
if $\geo G$ has local cut points then  $G$ 
splits over a virtually infinite cyclic subgroup unless $\geo G\simeq S^1$
and $G$ maps onto a Schwarz triangle group  with finite kernel.  
Therefore from now on we will assume that $\geo G$ has no local 
cut points.

A 1-dimensional, compact, metrizable, connected,
locally connected space $Z$ with no local cut points is a Menger
curve provided no point $z\in Z$ has a neighborhood which
embeds in the plane (see section \ref{prelim}).  Hence either $\geo G$ 
is a Menger curve or
some $\xi\in\geo G$ has a planar neighborhood $U$; therefore we
assume the latter holds.

\begin{lemma}
Let $\Ga\subset\geo G$ be a subset homeomorphic to a
finite graph.  Then $\Ga$ is a planar graph.
\end{lemma}

\proof
Since the action of $G$ on $\geo G$ is minimal, every $G$-orbit
intersects the planar neighborhood $U$, and so every point
of $\geo G$ has a planar neighborhood.
Because $\geo G$ has no local cut points, we have $\geo G\setminus \Ga
\neq \emptyset$.    So we can
find a hyperbolic element $g\in G$ whose fixed point set
$\{\eta_1,\eta_2\}\subset\geo G$ is disjoint from $\Ga$ (section \ref{prelim}).
Hence for sufficiently large $n$, $g^n(\Ga)$ is contained
in a planar neighborhood of $\eta_1$ or $\eta_2$.
\qed

We recall \cite{Claytor,Moise} that a compact, metrizable,
connected, locally connected space $X$ with no global cut points
is planar as long as no nonplanar graph embeds in $X$.
Therefore $\geo G$ is planar.  Finally, by \cite{Whyburn},
$\geo G$ is Sierpinski carpet.
\qed

\subsection{Groups with Sierpinski carpet boundary}

Let $M$ be a compact hyperbolic manifold with
nonempty totally geodesic boundary and let
$G\defeq\pi_1(M)$ be its fundamental group.
The universal cover $\tilde M$ of $M$ may be
identified with a closed convex subset of $\H^3$ which is bounded
by a countable disjoint collection ${\cal P}$ of  totally
geodesic planes.  Each $P\in{\cal P}$
bounds an open half-space disjoint from $\tilde M$.
$\tilde M$ is obtained from $\H^3$ by removing
each of these open half-spaces, and 
$\geo \tilde M\subset\geo\H^3$ is obtained from 
$\geo\H^3\simeq S^2$ by deleting the open disks 
corresponding to these half-spaces.  The closures
of these disks are disjoint since the distance 
between distinct elements of ${\cal P}$ is bounded 
away from zero.  As $\geo\tilde M$ has no interior
points in $S^2$, it is a Sierpinski carpet (see section \ref{prelim}).
Note that the peripheral circles of $\geo\tilde M$ 
are in 1-1 correspondence with elements of ${\cal P}$, 
and therefore the conjugacy classes
of $G$-stabilizers of peripheral circles are in 
1-1 correspondence with ${\cal P}/G$, the set
of boundary components of $M$.  The stabilizer of a 
peripheral circle is
the same as the stabilizer of the corresponding
element of ${\cal P}$, so these stabilizers
are quasi-convex in $G$.

The next theorem shows that similar conclusions hold
for any hyperbolic group whose boundary
is a Sierpinski carpet.

\begin{theorem}
Let $G$ be a hyperbolic group with boundary 
homeomorphic to the Sierpinski carpet $\S$.
Then

1. There are finitely many $G$-orbits of peripheral circles
in $\S$.

2. The stabilizer of each peripheral circle $C$ is a quasi-convex
 subgroup $G$ whose boundary is $C$.
 
\end{theorem}

\proof
We recall that $G$ acts cocompactly on the space 
$\trip G\defeq\{ (x,y,z)\in (\geo G)^3\mid \mbox{$x,\,y,\,z$ distinct}\}$.
Therefore if $C_k\subset\geo G$ is a sequence of
peripheral circles,  $(x_k,y_k,z_k)\in\trip G$ and
$\{x_k,y_k,z_k\}\subset C_k$, then after passing to 
a subsequence we may find a sequence $g_k\in G$,
$(x_\infty,y_\infty,z_\infty)\in\trip G$ so that
$(g_kx_k,g_ky_k,g_kz_k)$ converges to $(x_\infty,y_\infty,z_\infty)$.
But this means that $Diam(g_k(C_k))$ is bounded 
away from zero, so $g_k(C_k)$ belongs to a finite
collection of peripheral circles, and hence $g_k(C_k)$ is
eventually constant.
We conclude that there are only finitely many $G$-orbits of
peripheral circles, and the stabilizer of any $C\in{\cal C}$
acts cocompactly on the space of distinct triples in
$C$.  By \cite{Bowditch2}  $Stab(C)$ is a quasi-convex subgroup
of $G$, and $\geo Stab(C)=C$. From now on we will refer to stabilizers of 
peripheral circles as {\bf peripheral subgroups}. 
By \cite{Gabai, Casson, Tukia}
each peripheral subgroup is, modulo
a finite normal subgroup,  a cocompact Fuchsian
group in $Isom(\H^2)$.\qed

\subsection{Doubling Sierpinski carpet groups along peripheral 
subgroups}
\label{Doubling}

In this section we prove Theorem \ref{Doubletheorem}.

\medskip
Let $G$ be a hyperbolic group with $\geo G\simeq\S$,
and let $H_1,\ldots,H_k$ be a set of representatives
of conjugacy classes of peripheral subgroups of $G$.
We define a graph of groups ${\cal G}$ as follows.
The underlying graph has two vertices
and $k$ edges (no loops). Each vertex is labelled by a copy
 of $G$,
the $i^{th}$ edge is labelled by $H_i$, and the
edge homomorphisms $H_i\ra G$ are given by the inclusions.
We let $\hat G$ be the fundamental group of ${\cal G}$.

Next we construct a tree of spaces on which the group $\hat G$ acts in a natural way. 
Let $X_0$ be a finite Cayley 2-complex for $G$, and let $X_i$ be a finite 
Cayley 2-complex for the group $H_i$. The inclusion $H_i \hook G$ 
is induced by a cellular map $h_i : X_i \to X_0$ between the 
$2$-complexes.  Let 
$h:\cup X_i\ra X_0$ be the corresponding map from the disjoint union
of the $X_i$'s to $X_0$, and let $X$ denote the mapping
cylinder of $h$.

Let $DX$ be the double of $X$ along the collection of subcomplexes 
$X_i, i=1,...,k$. Consider now the universal cover $\dxt$ of $DX$ with the deck 
transformation 
group $\hat G$. Let $Y$ be the 1-skeleton of $\dxt$. 
The 1-skeletons of the subcomplexes $X_i, i=1,...,k$ lift to disjoint  
{\em edge subspaces} of $Y$. The {\em vertex subspaces} of $Y$ are 
closures of the connected components  to the complement to edge spaces. Each 
vertex space 
is the  cover of 1-skeleton of $X$. Let $T$ be the graph dual to the decomposition of 
$Y$ into vertex and edge subspaces: vertices $v$ of $T$ correspond to {\em vertex 
spaces} $Y_v\subset Y$, the edges $e$ correspond to the {\em edge subspaces} 
$Y_e \subset Y$. An edge $e$ is incident to a vertex $v$ if and only if $Y_e$
is contained in $Y_v$.  It is standard that the graph $T$ is actually a tree 
(compare \cite{ScottWall}). Let $V$ and $E$ denote the collections
of vertices and edges in $T$ respectively.
If $v\in T$ we let $E_v$ denote the collection
of edges containing $v$. 

 Let $\si:DX\ra DX$ be the natural involution of $DX$.
A map $\tau:Y\ra Y$ is a {\em reflection}
if it is a lift of $\si$ and it fixes some point;
each reflection fixes some edge space in $Y$,
and each edge space $Y_e$ is the fixed point set of 
precisely one reflection $r_e$. 
Let $\Ga$ be the group generated by the reflections
in $Y$. The group $\Ga$ is normalized by $\hat G$ since conjugation
of a reflection by an element of $\hat G$ yields
another reflection; likewise $\hat G$ is normalized by $\Ga$.
Let $v\in T$ be any vertex.  Then
$\Ga$ is the free product of order two subgroups
of the form $\langle r_e\rangle$ where $e\in E_v$.
The vertex space $Y_v$ is a fundamental domain for the
action of $\Ga$ on $Y$. The group $\Ga$ preserves the tree
structure of $Y$, so we have an induced action
of $\Ga$ on $T$ by tree automorphisms, each 
reflection $r_e$ acting on $T$ as an inversion
of the edge $e$.  The action of $\Ga$ on 
$T$ naturally induces an action of $\Ga$ on
$\geo T$. The space $Y$ is a connected graph, and we give it
the 
natural path-metric where each edge in $Y$ has unit length.  

\begin{lemma}
\label{edgevertexprops}
1. The space $Y$ is Gromov-hyperbolic. 

2. Edge and vertex spaces are all $K$-quasi-convex 
in $Y$ for some $K$. 

3. There is a function $C(R)$ such that 
for every $R$, the intersection of $R$-neighborhoods 
of any two distinct vertex or edge spaces has diameter at most $C(R)$ 
unless the spaces are incident.
\end{lemma}

\proof The space $Y$ is quasi-isometric to Cayley graph of $\hat G$. The group 
$\hat G$ is Gromov-hyperbolic by \cite{BF2, BF3}. The assertions 2 and 3 follow from \cite{Mitra} and \cite{Swarup2}. 
\qed

\medskip
We have a coarse Lipschitz projection $p:Y\ra T$
which maps $(Y_v-\cup_{e\in E_v}Y_e)$ to $v$ for each $v\in V$,
and maps each edge space to the midpoint of the corresponding
edge of $T$.  
If $\ga:[0,\infty)\ra \geo Y$ is a 
unit speed geodesic ray, then $p\circ\ga$ is a coarse Lipschitz
path with the bounded backtracking property\footnote{A
map $c:[0,\infty)\to T$ has the {\bf bounded backtracking
property} if for every $r\in (0,\infty)$ there is an $r'\in (0,\infty)$ such 
that if $t_1<t_2$, and $d(c(t_1),c(t_2))>r'$, then $d(c(t),c(t_1))>r$
for every $t>t_2$.}
 by 
the quasi-convexity of vertex/edge spaces. Hence   $p\circ\ga$ 
maps into a finite tube around a geodesic ray $\tau$ in $T$. 
If $p\circ \ga$ is unbounded in $T$, then the equivalence class
of the ray $\tau$ is 
uniquely determined by $\ga$ 
and we label $\ga$ with the associated boundary point
$[\tau]\in\geo T$.
By the quasi-convexity of edge spaces, if $\ga$ hits an 
edge space for an unbounded sequence of times, then 
it remains in a quasi-convex neighborhood of the edge space.
In this case, we know that $\ga$ eventually remains in
a bounded neighborhood of a unique edge space by property 3 in Lemma 
\ref{edgevertexprops}, and we label $\ga$ with this edge.  
If neither of the above two cases occurs, then for each edge $e$
of the tree, we know that $\ga$ eventually lies in
one of the two components of the complement of the edge
space $Y_e$, and we label the edge with 
an arrow pointing in the direction of the corresponding
subtree of $T$.  There must be some (and at most one) 
vertex $v\in T$  such that all edges emanating from $v$
have arrows pointing toward $v$; otherwise we could follow
arrows and leave any bounded set.   There must be an unbounded
sequence of times $t_k$ such that $\ga(t_k)$ lies in the
vertex space $Y_v$ (by the construction of the edge labelling);
by quasi-convexity of $Y_v$, this means that $\ga$ eventually
lies in the $R$-neighborhood of $Y_v$; in this case we label $\ga$ by $v$.
Equivalent geodesic rays are given the same label.
We get a labelling map $Label:Y\ra (T\cup\geo T)$
which is clearly $\Ga$-equivariant.

We now examine the topology of $\geo Y$. This space is metrizable and we 
fix a metric $d$ on $\geo Y$; in what follows we will implicitly 
use $d$ when discussing metric properties of $\geo Y$. 
Recall that each vertex space $Y_v$ is 
quasi-isometric to $G\simeq\tilde X$; since
by Lemma \ref{edgevertexprops} every subspace $Y_v$ is 
quasi-convex in $Y$, we conclude that $\geo Y_v\subset\geo Y$
is a Sierpinski carpet.  
Similarly, the peripheral circles of the Sierpinski carpet
$\geo Y_v$ are in 1-1 correspondence with
the boundaries of edge spaces $Y_e\subset Y_v$.  

By the visibility property of the uniformly quasi-convex 
edge spaces, there is at most one boundary point of $\geo Y$
labelled by any $\xi\in\geo T$.  For each edge $e$ in $T$,
the set of points in $\geo Y$ labelled by $e$ is the
ideal boundary of the edge space $Y_e$, i.e. a circle.
For each vertex $v\in T$, the set of points labelled
by $v$ is 
$$
\geo Y_v- \cup_{e\in E_v}\geo Y_e
$$
i.e. the Sierpinski carpet $\geo Y_v$ minus the union of its peripheral
circles.  

Our next goal is to describe the topology of $\geo Y$ using the tree $T$. 
Choose  $v\in T$. Every edge $e$ of
$T$ separates $T$ into two subtrees, and we let
$T_{v,e}\subset T$ be the subtree disjoint from $v$.
We define the {\bf outward subset}, $Out_{v,e}$, for a pair  $(v,\,e)\in V\times E$ 
to be the
collection of points of $\geo Y$ labelled by
elements of $T_{v,e}\cup \geo T_{v,e}$. 
The visibility property of $Y$ implies that for a fixed
$v\in T$ and any $\eps>0$
there are only finitely many edges $e\subset T$ so that
the diameter of $Out_{v,e}$ exceeds $\eps$.
Outward subsets of $\geo Y$ are open since a geodesic ray
$\ga$ with $\geo\ga\in Out_{v,e}$ will eventually leave
any tubular neighborhood of the edge space $Y_e$,
and so nearby boundary points correspond to rays
which eventually lie in the same component of the
complement of $Y_e$ in $Y$.  It follows that 
if $\xi\in\geo T$, and $e_k$ is the sequence of 
edges occurring in the ray $\ol{v\xi}$, then
the sequence of outward sets $Out_{v,e_k}$ is a 
nested basis for the topology of $\geo Y$ at
the point labelled by $\xi$. The closure of 
$Out_{v,e}$ is $Out_{v,e}\cup \geo Y_e$ because
the complement to $Out_{v,e}\cup \geo Y_e$
is $Out_{w,e}$ where $w$ is the endpoint of 
$e$ furthest from $v$ (obviously $\geo Y_e\subset
\ol{Out_{v,e}}$).

\begin{lemma}
\label{toplem}
Suppose $\xi_k\in\geo Y$ converges to $\xi_\infty\in\geo Y$.
Then one of the following holds.

1. $\xi_\infty$ is labelled by a boundary point 
$Label(\xi_\infty)\in\geo T$. In this case 
$Label(\xi_k)$ converges to $Label(\xi_\infty)$ in the
compact space $T\cup\geo T$.

2.  $\xi_\infty$ is labelled by a vertex $v\in T$.
In this case,
 for any subset ${\cal E}\subseteq E_v$
containing all but finitely many elements of $E_v$, 
the sequence $\xi_k$
eventually lies in $$\geo Y_v
 \cup(\cup_{e\in{\cal E}}Out_{v,e}).$$
 
3. $\xi_\infty$ is labelled by an edge $e_0$.
In this case, if $v,\,w$ are the endpoints of $e_0$
 then for any subset ${\cal E}\subseteq E_v$
containing all but finitely many elements of $E_v$,
and  any subset ${\cal F}\subseteq E_w$
containing all but finitely many elements of $E_w$,
the sequence $\xi_k$ eventually lies in
$$\geo Y_{v}\cup\geo Y_{w}\cup (\cup_{e\in{\cal E}}Out_{v,e})
\cup(\cup_{e\in{\cal F}}Out_{w,e}).$$ 

\end{lemma} 

\proof  {\em Case 1.} If $v$ is any arbitrary vertex of  $T$,
and $e_1,\,e_2,\ldots$ is the sequence of edges 
comprising the geodesic ray $\ol{v\xi_\infty}\subset T$,
then $Out_{v,e_j}\subset\geo Y$ is a neighborhood basis
for $\xi_\infty$.  Therefore $Label(\xi_k)$ converges to
$Label(\xi_\infty)$ by the definition of the topology on
$T\cup\geo T$. 

{\em Case 2.} If this weren't the case, then a subsequence
of $\xi_k$ would converge to an element of $\ol{Out_{v,e}}
=Out_{v,e}\cup \geo Y_e$
for some $e\notin {\cal E}$.  This contradicts the fact
that $\xi_\infty$ is labelled by $v$.

{\em Case 3.} Similar to case 2.

\qed

\begin{proposition}
$\geo\hat G$ is homeomorphic to $S^2$.  
\end{proposition}

\proof
Let $G'$ be the fundamental group of a compact
hyperbolic $3$-manifold $M$ with nonempty totally
geodesic boundary.  Recall (see section \ref{prelim})
that $\geo G'$ is a Sierpinski carpet.  Using 
the notation developed
above (decorated with ``primes''),  $\hat G'$ is the 
fundamental group of the
double of $M$, so $\geo \hat G'$ is homeomorphic
to $S^2$.  We will construct a homeomorphism
between $\geo\hat G'$ and $\geo \hat G$.

Choose vertices $v\in T$ and $v'\in T'$, 
and a bijection $E_v\ra E_{v'}$.
This induces an isomorphism between Coxeter
groups $\Ga\ra\Ga'$, which we will use to identify
$\Ga$ with $\Ga'$.  There is a unique
$\Ga$-equivariant isomorphism $T\cup\geo T
\ra T'\cup\geo T'$ which induces the given
bijection $E_v\ra E_{v'}$; we will use primes
to denote corresponding edges and vertices.  Choose
an enumeration $v=v_1,v_2,\ldots$ of vertices of 
$T$ so that $d(v_k,\cup_{j<k}v_j)=1$.  Choose
a homeomorphism $f_1:\geo Y_v\ra\geo Y'_{v'}$.  Using
reflections from $\Ga$ we inductively extend $f_1$ to a 
homeomorphism $f_k:\cup_{i=1}^k\geo Y_{v_i}\ra
\cup_{i=1}^k\geo Y'_{v_i'}$ for each $k$, so that 
the resulting map $f_\infty:\cup_{i=1}^\infty\geo Y_{v_i}\ra
\cup_{i=1}^\infty\geo Y'_{v_i'}$ is $\Ga$-equivariant.
By construction, $f_\infty$ is compatible with 
label maps, i.e. the following diagram commutes:

$$\begin{array}{ccc}
\cup_{i=1}^\infty\geo Y_{v_i} & \stackrel{f_{\infty}}{\lra} &  
\cup_{i=1}^\infty\geo Y_{v_i'}'  \\
Label\Big\downarrow & ~ &  Label\Big\downarrow\\
T\cup \geo T  & \stackrel{id}{\lra}  &  T\cup \geo T \\
\end{array}$$

We claim that $f_\infty$ extends continuously to a
homeomorphism $f:\geo Y\ra\geo {Y'}$.
In view of the naturality of our construction
it is enough to show that $f_\infty$ extends
to a continuous map $f:\geo Y\ra\geo Y'\simeq\geo \hat{G}'\simeq S^2$,
since the inverse map may be produced by 
exchanging the roles of $G$ and $G'$.  Pick
a sequence $\xi_k\in\geo Y$ which converges to 
some $\xi\in\geo Y$.  We will show that 
$f_\infty(\xi_k)$ converges.

{\em Case 1: $\xi$ is labelled by some $\eta\in\geo T$.}
In this case there is a unique $\xi'\in \geo {Y'}$
which is labelled by $\eta'\in\geo T'$.
We know that if $e_i$ (resp $e_i'$) is the sequence
of edges of the ray $\ol{v\eta}$ (resp $\ol{v'\eta'}$), then the
outward sets $Out_{v,e_i}$ (resp. $Out_{v',e_i'}$)
form a basis for the topology of $\geodx$
(resp. $\geo {Y'}$) at $\xi$ (resp. $\xi'$).
Since $f_\infty$ maps $Out_{v,e_i}\cap 
\cup_{i=1}^\infty\geo Y_{v_i}$ to 
$Out_{v',e_i'}\cap 
\cup_{i=1}^\infty\geo Y_{v_i'}'$, the sequence 
$f_\infty(\xi_k)$ converges to $\xi'$.

{\em Case 2: $\xi$ is labelled by a vertex $v\in T$.}
For each $k$ either $\xi_k\in\geo Y_v$ or $\xi_k\in
Out_{v,e_k}$ for a unique $e_k\in Edge_v$. By Lemma
\ref{toplem}, in the latter case $Diam(Out_{v,e_k})\ra 0$
as $k\ra\infty$.  Construct a sequence 
$\zeta_k\in\geo Y_v$ so that $\zeta_k=\xi_k$
when $\xi_k\in\geo Y_v$, and 
$\zeta_k\in\geo Y_{e_k}=\ol{Out_{v,e_k}}\cap\geo Y_v$
otherwise.  Note that $\lim_{k\ra\infty}\zeta_k=\xi$
since $Diam(Out_{v,e_k})\ra 0$.  The sequence
 $f_\infty(\zeta_k)$ converges to $f_\infty(\xi)$
 since $f\restr_{\geo Y_v}$ is continuous.  Observe that
 $d(f_\infty(\zeta_k),f_\infty(\xi_k))$ is
 zero when $\xi_k\in\geo Y_v$ and 
 is at most $Diam(Out_{v',e_k'}')$ otherwise. 
 Since each $e_k$ occurs only finitely often,
 $Diam(Out_{v',e_k'}')\ra 0$ so 
 $$\lim_{k\ra\infty}f_\infty(\xi_k)=\lim_{k\ra\infty}
 f_\infty(\zeta_k)=f_\infty(\xi).$$

{\em Case 3: $\xi$ is labelled by an edge $e_0\in T$.}
We leave this case to the reader, as it is similar to 
case 2.

\qed

\begin{corollary}
\label{pd3pair}
If $G$ is torsion-free, then so is $\hat G$, 
and in this case $\hat G$ is a $3$-dimensional 
Poincar\'e duality group
by  \cite{BM}, \cite{Best}.  By \cite{DD}, if one splits 
a $PD(n)$ group over a 
$PD(n-1)$ subgroup, then the vertex groups (together
with the incident edge subgroups) define relative
$PD(n)$ pairs; therefore  $(G;H_1,\ldots,H_k)$
is a relative Poincar\'e duality pair.  
In particular $\chi(G)=\frac{1}{2}\sum_i\chi(H_i)<0$.
\end{corollary}

\begin{corollary}
\label{cannonimpliesconj}
Let $G$ be a torsion-free hyperbolic group with Sierpinski
carpet boundary.  Suppose either

A. Cannon's conjecture is true

or

B. Every $3$-dimensional Poincar\'e duality group with a nontrivial splitting is the  
fundamental group of a closed $3$-manifold.

\no
Then $G$ is the fundamental group of a compact hyperbolic
$3$-manifold with totally geodesic boundary.
\end{corollary}

\proof
Let $H_1,\ldots,H_k$, $\hat G$, $\Ga$, 
be as in the  first part of this section.  If A holds,
then $\hat G$ is the fundamental group of a closed
hyperbolic $3$-manifold  $M$.  Since $\hat G$
splits nontrivially by its very definition, if B holds
then $\hat G=\pi_1(M)$ where $M$ is a closed irreducible
$3$-manifold.  $M$ is Haken since its fundamental group
splits, and so Thurston's uniformization theorem implies
that $M$ admits a hyperbolic structure.  In either
case we have $\hat G$ acting on $\H^3$ discretely,
cocompactly, and isometrically.

The reflection group $\Ga$ acts on $\hat G$ by conjugation,
with each reflection centralizing a unique quasi-convex
edge subgroup of $\hat G$.
By Mostow rigidity, $\Ga$ acts isometrically on the universal
cover of $M$ normalizing the action  $\hat G\acts\H^3$.  
$G\subset\hat G$ is a quasi-convex subgroup, and so 
it acts on $\H^3$ as a convex cocompact subgroup.   
The limit set of $G$ in $\geo\H^3$ is a Sierpinski carpet,
and because every peripheral subgroup of $G$ is 
centralized by a unique reflection in $\Ga\subset Isom(\H^3)$, 
the peripheral circles are fixed by reflections in 
$\Ga$.  Thus each peripheral circle of the limit set of $G$
is a round circle, and so the convex hull of the limit set
is a convex subset bounded by disjoint totally geodesic
hyperbolic planes.  It follows that $G$ is the fundamental
group of a compact hyperbolic manifold with totally geodesic
boundary.

\qed

\subsection{Examples}

We now use Theorems 1 and \ref{Doubletheorem} to see that some
classes of hyperbolic groups have Menger curve boundary.

We first remark that a torsion-free hyperbolic group with Sierpinski
carpet boundary has negative Euler characteristic by 
Corollary  \ref{pd3pair}.  So if $G$ is a torsion-free hyperbolic group with 1-dimensional boundary,
$G$ doesn't split over a trivial or cyclic group, and $\chi(G)\geq 0$,
then $\geo G$ is a Menger curve.

\begin{theorem}
Let $G$ be a torsion-free 2-dimensional hyperbolic group that does not 
split over trivial and cyclic subgroups and which fits into a short exact 
sequence:
$$
1\to F \to G \to \Z \to 1
$$
where $F$ is finitely generated. Then $\geo G$ is the Menger curve. 
\end{theorem}
\proof In view of Theorem 1, it is enough to show that $\geo G$ is not 
Sierpinski carpet. Suppose it is. Note that if $F$ admits a finite
Eilenberg-Maclane space, then it is easy to see that $\chi(G)=\chi(F)\chi(\Z)=0$,
so $\geo G$ cannot be a Sierpinski carpet by the remark above. 
However there are examples such that $F$ is not a finitely presentable group  
(see \cite{Rips}). We now consider the general case.  
Then $(G; H_1,\ldots, H_k)$ is a {\em relative Poincare duality pair}. 
Let $K_0$ be a finite
Eilenberg-Maclane space for the group $G$, let $D$ be a disjoint
union of finite Eilenberg-Maclane spaces for the groups $H_1,\ldots,H_k$,
and let $K$ be the mapping cylinder for a map $D\ra K_0$
which induces the given maps $H_i\hook G$.  We view $D$ as a subcomplex
of $K$.  Consider the finite cyclic coverings
$$
(K_n, D_n)\to (K, D)
$$
which are induced by the homomorphisms $G\to \Z \to \Z_n$. Then each pair 
$(K_n, D_n)$  again satisfies relative Poincare duality in dimension 3, so
$$
H^*(K_n, D_n; \Z/2) \cong \tilde{H}^{3-*}(K_n; \Z/2)
$$  
We will use the notation $b_j(L)$ to denote the rank of $H_j(L, \Z/2)$. 
Thus $\lim_{n\to \infty} b_1(D_n)= \infty$ and 
$b_1(K_n)\le b_1(F) +1< \infty$.  Consider the exact sequence 
of the pair $(K_n, D_n)$:
$$
... \to H^1(K_n;\Z/2) \to H^1(D_n;\Z/2)\to H^2(K_n, D_n;\Z/2) \to ...
$$
Since $b_1(K_n)$ is bounded by $b_1(F)+1$ and 
$\lim_{n\to \infty} b_1(D_n)= \infty$, it follows that 
$\lim_{n\to \infty} Dim(H^2(K_n, D_n;\Z/2))=\infty$. 
This contradicts the fact that  $H^2(K_n, D_n;\Z/2) \cong H_1(K_n;\Z/2)$. \qed

\medskip
Now let $F$ be a finitely generated free group and $\phi: F\to F$ be a {\bf
hyperbolic} automorphism (see \cite{BF2} for the definition). 
Consider the extension
$$
1\to F \to G \to \Z \to 1
$$
induced by $\phi$. The group $G$ is hyperbolic by \cite{BF2}. The cohomological dimension of $G$ is 2 by the Mayer-Vietoris sequence, thus the boundary of $G$ is 1-dimensional by \cite{BM}.

\begin{corollary}
$\geo G$ is the Menger curve. 
\end{corollary}
\proof We will show that the group $G$ does not split over a cyclic (possibly trivial) subgroup. Suppose that it does. Then we have the corresponding action of $G$ on a minimal simplicial tree $T$ with cyclic edge stabilizers. Consider 
the restriction of this action on the subgroup $F$. Let $T'\subset T$ be the minimal $F$-invariant subtree, then $T'$ is $\Z$-invariant (since $\Z$ normalizes $F$), thus $T'=T$. By Grushko's theorem (in the case of trivial 
edge stabilizers) and the generalized accessibility theorem \cite{BF1} (in the case of infinite cyclic stabilizers), the quotient  $T/F$ is a finite graph $\Ga$. The action of $\Z=\<z\>$ projects to action on $\Ga$, after taking a 
finite iteration of  $\phi$ (if necessary) we may assume that $z$ acts 
trivially on $\Ga$. Since $G$ does not contain $\Z^2$-subgroups, the edge stabilizers for the action of $F$ on $T$ must be trivial. 
Thus we get a free product decomposition of $F$ so that each factor is invariant under some iterate of $z$. This contradicts the assumption that the corresponding automorphism $\phi: F\to F$ is hyperbolic. 
\qed

\begin{theorem}
Let ${\cal G}$ be a finite graph of groups.  Suppose

1. Each vertex group is a torsion-free hyperbolic group
whose boundary is either a Menger curve or a Sierpinski carpet;
and at least one vertex group has Menger curve boundary.

2. Each edge group is a finitely generated free group of rank at
least 2, and includes as a quasi-convex subgroup of each of the 
corresponding vertex groups.

3. If $T$ is the Bass-Serre tree for ${\cal G}$, and $e_1,\,e_2\subset T$
are two edges emanating from the same vertex $v\in T$, then 
their stabilizers intersect trivially.

\no
Then the fundamental group $G$ of ${\cal G}$ is a hyperbolic group with
Menger curve boundary.
\end{theorem}

\proof
Conditions 2 and 3 imply that $G$  is hyperbolic by \cite{BF2},
and vertex groups are quasi-convex subgroup of $G$ by \cite{Mitra, Swarup2}.
$G$ is torsion-free since all vertex groups are torsion-free.
$G$ has cohomological dimension 2 by the Mayer-Vietoris sequence, so $\geo G$
has dimension $1$ by \cite{BM}.  

We claim that $G$ does not
split over trivial or infinite cyclic groups.  To see this,
let $T$ be the Bass-Serre tree of ${\cal G}$, and let $S$
be the Bass-Serre tree of a splitting of $G$ over  trivial
and/or cyclic groups.  Consider two adjacent vertices $v_1,\,v_2\in T$,
 let $G_{v_i}\subset G$ be their stabilizers, and let $G_e$
 be the stabilizer of the edge joining them.  Since
$G_{v_i}$ does not split over trivial or cyclic subgroups \cite{Bowditch2},
$G_{v_i}$ has a nonempty fixed point set in $S$.  If $s_i\in S$
is fixed by $G_{v_i}$, then the segment joining $s_1$ to $s_2$
will be fixed by $G_e$.   Since $G_e$ is free of rank at least $2$,
we see that $s_1=s_2$.  Therefore by induction we find that
$G$ has a global fixed point in $S$, which is a contradiction.

If the stabilizer of $v\in T$ has Menger curve boundary,
then by the quasi-convexity of $G_v$ in $G$, the Menger
curve embeds in $\geo G$.  This shows that $\geo G$ cannot
be homeomorphic to $S^1$ or the Sierpinski carpet.  By 
Theorem \ref{s1mengersierpinski} $\geo G$ is a Menger curve.
\qed

\subsection{Topologically rigid groups}
\label{toprigid}

In this section we will construct some examples of 
topologically rigid groups.  Before proceeding, we first note
 a consequence of Theorem \ref{s1mengersierpinski}.

 \begin{corollary} Let $G$ be a nonelementary hyperbolic group with
 $Dim(\geo G)\leq 1$.  Then $G$ is not topologically rigid.
 \end{corollary} 
 
 \no
 We will sketch a proof of the corollary, and leave the details
 to the reader.
 
 {\em Case I: $G$ has more than one end.}  Then $G$ splits as an
 amalgamated product or HNN extension over a finite
 group. Let $G\acts T$ be the action of $G$ on the Bass-Serre
 tree associated to such a splitting, so there is only
 one edge orbit in $T$.  Following along the same lines as in section
 \ref{Doubling}, we construct a tree of spaces $X$, with vertex
 and edge spaces corresponding to vertices and edges in $T$.
 For each vertex $v\in T$, the vertex space  $X_v\subset X$
 is quasi-convex in $X$ and as in section \ref{Doubling}
 we may label points in $\geo X$ with elements of
 $T\cup\geo T$.  The outward sets (see section \ref{Doubling})
 are open and closed in $\geo X$.  If $e_1$ and $e_2$ are incident
 to a vertex $v$ then they lie in the same $G_v$-orbit (since $G/T$
 has only one edge).   $Out_{v,e_1}$ and $Out_{v,e_2}$ are
 disjoint and homeomorphic, so we may define a homeomorphism
 of $\geo X$ by swapping them while holding everything else
 fixed.  This construction yields a continuum of homeomorphisms
 of $\geo X$, so $G\ra Homeo(\geo X)$ cannot be surjective.
 
 {\em Case II: $G$ is 1-ended.}  If $\geo G$ is homeomorphic
 to $S^1$, the Sierpinski carpet, or the Menger curve then $G$ cannot
 be topologically rigid since each of these spaces has uncountable homeomorphism
 group.  Therefore by Theorem \ref{s1mengersierpinski} we may assume
 that $G$ splits as an amalgamated free product or HNN extension
 over a virtually cyclic group.  Let $G\acts T$
 be the action of $G$ on the Bass-Serre tree associated with such
 a splitting.
 If $e$ is an edge in $T$, $e=\ol{v_1v_2}$, then 
 $Out_{v_1,e}-\geo X_e$ and $Out_{v_2,e}-\geo X_e$ are open
 and closed in $\geo X-\geo X_e$, and are preserved by $G_e$.
 Take a $g\in G_e$ so that $\geo G$ fixes both points in 
 $\geo G_e$, and define a homeomorphism $f:\geo X\ra\geo X$
 by $f\restr_{Out_{v_1,e}}=\geo g\restr_{Out_{v_1,e}}$
 and $f\restr_{Out_{v_2,e}}=id\restr_{Out_{v_2,e}}$.  This 
 type of construction will give a continuum of homeomorphisms
 of $\geo X$, so again $G\ra Homeo(\geo X)$ cannot be surjective.
 \qed

 Our construction of topologically rigid groups is based on the
 idea (realized precisely in Proposition \ref{stabismobius})
 that a homeomorphism of $S^2$ must be a M\"obius transformation
 provided it preserves a sufficiently rich family of round circles.
 We begin with an analogous statement for homeomorphisms of $S^1$. 
 
\bigskip
\no
{\bf Line configurations in $\H^2$.} 
Let ${\cal L}$ be a locally finite collection of 
geodesics in $\H^2$ so that  the complementary regions of $\cup_{L\in{\cal L}}L$ are bounded, and we assume that there is a cocompact lattice   
 $\Ga\subset Isom(\H^2)$ stabilizing ${\cal L}$.
  Let $\db\H^2$ be the space of unordered distinct pairs
in $\geo \H^2$, and let $\geo{\cal L}$ be the collection of
pairs of endpoints $\geo L$ for $L\in {\cal L}$,
$\geo{\cal L}\defeq \{\geo L\mid L\in{\cal L}\}\subset\db\H^2$.
Note that if $L_1,\,L_2\in{\cal L}$ and $\geo L_1\cap\geo L_2\neq\emptyset$
then $L_1=L_2$.  Let $Stab(\geo{\cal L})\subset Homeo(\geo\H^2)$
be the group of homeomorphisms of $\geo\H^2$ which preserve
$\geo{\cal L}\subset\dub\H^2$.

\begin{lemma}
\label{lines}

1. If  $L_1,\,L_2\in{\cal L}$
have nonempty intersection and $g\in Stab(\geo{\cal L})$ fixes $\geo L_1
\cup \geo L_2$ pointwise then $g=id$.  

2. $\{ \geo\ga\mid \ga\in \Ga\}
\subset Homeo(\geo\H^2)$  is a finite index
subgroup of $Stab(\geo {\cal L})$. 
\end{lemma}

\proof  Our arguments essentially follow \cite[Proof of Theorem 2.7]{Cass}. 
We will identify the space of geodesics in $\H^2$ with $\dub\H^2$.

(1) Suppose $L_1,\,L_2\in{\cal L}$ and $g\in Stab(\geo{\cal L})$
fixes $\geo L_1\cup\geo L_2$ pointwise.
If $\si_1,\,\si_2$ are the connected components of 
$\geo\H^2-\geo L_1$, then $g(\si_i)=\si_i$ since
$|\geo L_2\cap \si_i|=1$ and $\geo L_2$ is fixed by $g$.   Observe that
$\Si_i\defeq\{\geo L\cap\si_i\mid \mbox{$L\in {\cal L}$ and $|L\cap L_1|=1$}\}
\subset \si_i$ is a discrete subset of $\si_i$ with the order type
(with respect to the ordering on $\si_i\simeq\R$) of the integers,
and $g(\Si_i)=\Si_i$.   But $g$ fixes the point $\geo L_2\cap \si_i\in\Si_i$
and is orientation preserving, so $g\restr_{\Si_i}=id_{\Si_i}$.
Therefore $g$ fixes $\geo L$ for every $L\in{\cal L}$ with $L\cap L_1\neq\emptyset$.
The incidence graph of ${\cal L}$ is connected, so we may apply
this argument inductively to see that $g$ fixes $\geo L$ for 
every $L\in {\cal L}$.  The set $\cup_{L\in{\cal L}}\geo L$
is dense in $\geo \H^2$, so $g=id$.  This proves the first assertion
of the lemma.

(2) We now show that every sequence $g_k\in Stab(\geo{\cal L})$ has
a subsequence which is constant modulo $\Ga$, which proves
that $[Stab(\geo {\cal L}):\Ga]<\infty$. Pick $L_1,\,L_2\in{\cal L}$
such that $L_1$ intersects $L_2$ in a point $p$.  For each
$k$ let $g_{k*}L_i\in{\cal L}$ be the unique line with
$\geo (g_{k*}L_i)=g_k(\geo L_i)$.  Then 
$(g_{k*}L_1)\cap (g_{k*}L_2)=p_k$ for some $p_k\in\H^2$,
and we may choose a sequence $\ga_k\in\Ga$ such that
$\sup d(\ga_k(p_k),p)=R<\infty$.  Then the lines $(\ga_k\circ g_k)_*L_i$
lie in the  finite set $\{ L\in{\cal L}\mid {\cal L}\cap \ol{B(p,R)}\neq
\emptyset\}$,
so after passing to a subsequence we may assume that 
$(\ga_k\circ g_k)\restr_{\geo L_i}$ independent of $k$ for $i=1,\,2$.
By the previous paragraph the sequence $\ga_k\circ g_k\in Homeo(\geo \H^2)$
is constant.
\qed

\medskip
\no
{\bf Plane configurations in $\H^3$.}
Below we prove an analog of Lemma \ref{lines} for a  collection ${\cal H}$ of totally 
geodesic hyperplanes in $\H^3$.

Let ${\cal H}$ be a locally finite collection of totally geodesic 
planes in $\H^3$, with stabilizer $G\defeq\{ g\in Isom(\H^3)\mid
\mbox{$g(H)\in{\cal H}$ for every $H\in{\cal H}$}\}$.  Let $\geo{\cal H}
\defeq \{\geo H\mid H\in{\cal H}\}$.
We assume
that ${\cal H}$ satisfies the conditions:

1. $G$ is a cocompact lattice in $Isom(\H^3)$.

2. The complementary regions of $\cup_{H\in{\cal H}}H$ are bounded.

3. If $H\in {\cal H}$, then reflection in $H$ does not preserve the collection
${\cal H}$.

\no
Such examples will be constructed later in this section. 

\smallskip
The local finiteness of ${\cal H}$ implies that there are finitely
many $G$-orbits in ${\cal H}$, and that the stabilizer of each
$H\in{\cal H}$ acts cocompactly on $H$.

\begin{definition}
We will say that three circles $\geo H_1,\,\geo H_2,\, \geo H_3$,
where $H_i\in{\cal H}$, are in {\bf standard position} if 
the three planes $H_i$ intersect transversely in a single
point $x\in \H^3$.
\end{definition}

Note that if the circles $\geo H_1,\,\geo H_2,\, \geo H_3$ are in 
standard position and $C_1, C_2, C_3$ is another unordered triple 
of circles which 
bound elements of ${\cal H}$, then $C_1, C_2, C_3$ is in standard 
position if and only if there is a homeomorphism 
$f:  \geo H_1 \cup\geo H_2 \cup \geo H_3 
\to C_1 \cup C_2 \cup C_3$ which carries elements of ${\cal H}$ to 
elements of ${\cal H}$.

Let $Stand$ denote the collection of unordered triples of circles
in standard position.  We will say that two elements of $Stand$ are 
{\bf incident} if they have exactly two circles in common.

\begin{lemma}
\label{circlescut}
1. The incidence graph of $Stand$ is connected.

2. If $\ga\subset\geo\H^3$ is homeomorphic to $S^1$, then
either $\ga=\geo H$ for some $H\in {\cal H}$, or there is
an $H\in {\cal H}$ so that $\geo H$ intersects both components
of $\geo \H^3-\ga$.
\end{lemma}

\proof
The union $\cup_{H\in{\cal H}}H$ determines a 
polygonal subcomplex in $\H^3$ with connected 1-skeleton.
Therefore the assertion 1 follows.

To prove the assertion 2, let $U$ and $U'$ denote the connected components
of $\geo \H^3-\ga$.  We may find $H,\,H'\in{\cal H}$ so that
$\geo H\subset U, \geo H'\subset U'$.  Since the incidence graph for
${\cal H}$ is connected we can find a chain of planes $H_0=H, H_1,..., H_n=H'$ 
in ${\cal H}$ so that consecutive planes intersect each other. 
We see that either $\ga=\geo H_j$
for some $H_j$ in this sequence or for some $H_j$
the circle $\geo H_j$ intersects both $U$ and $U'$.
\qed

\begin{proposition}
\label{stabismobius}
Let $Stab(\geo{\cal H})$ be the group of homeomorphisms of $\geo\H^3$
which preserve $\geo{\cal H}$, $Stab(\geo{\cal H})\defeq
\{ g\in Homeo(\geo \H^3)\mid 
\mbox{$g(\geo H)\in\geo {\cal H}$ for all $H\in {\cal H}$}\}$. Then 
$Stab(\geo{\cal H})=\{\geo g\mid g\in G\}$.
\end{proposition}
\proof 
Suppose $\{\geo H_1,\,\geo H_2,\,\geo H_3\}\in Stand$, 
$f\in Stab(\geo {\cal H})$, and $f(\geo H_i)=\geo H_i$
for $1\leq i\leq 3$.  Then for $1\leq i\leq 3$ we may
consider the collection ${\cal L}_i$ of geodesics in 
$H_i$ of the form $H_i\cap H$ for $H\in {\cal H}-H_i$.
Part 1 of Lemma \ref{lines} then implies that $f\restr_{\geo H_i}=id_{\geo H_i}$.

Now suppose $\{\geo H_1,\geo H_2,\geo H_3\},\,\{\geo H_1,\geo H_2,\geo H_4\}\in Stand$
are incident, $f\in Stab(\geo {\cal H})$, and $f\restr_{\geo H_i}
=id_{\geo H_i}$ for $1\leq i\leq 3$.  Then $f(\geo H_4)=\geo H_4$
since $H_4$ is the unique element of ${\cal H}$ whose boundary
contains the 4 element set $\geo H_4\cap(\geo H_1\cup\geo H_2)$.
Therefore by the previous paragraph we have 
$f\restr_{\geo H_4}=id_{H_4}$.  Since the incidence graph of
$Stand$ is connected we see by induction that $f\restr_{\geo H}=id_{\geo H}$
for all $H\in {\cal H}$, and this forces $f=id_{\geo \H^3}$.

Reasoning as in Lemma \ref{lines} we conclude that 
$[Stab(\geo{\cal H}):G]<\infty$.  

Let $G'\subset G$ be a finite index normal subgroup
of $Stab(\geo{\cal H})$.  Each $f\in Stab(\geo {\cal H})$ 
normalizes the action $G'\acts\geo\H^3$, so by Mostow rigidity 
each $f$ is a M\"obius transformation.
Therefore for every $f\in Stab(\geo {\cal H})$ we have
$f=\geo g$ for some $g\in G$.
\qed

\medskip 
\no
{\bf Constructing topologically rigid groups.}
 Let $G'\subset G$ be a finite index torsion-free subgroup
 of $G$.   Let $\{H_1,\ldots,H_k\}$
be a set of representatives of the $G'$-orbits in ${\cal H}$,
and let $G_i\defeq Stab(H_i)$.  For any $1\leq i\leq k$, the
set of geodesics $\{ H\cap H_i\mid H\in {\cal H}-H_i,\,H\cap H_i\neq\emptyset
\}\subset H_i$ is
finite modulo the action of $G_i$.  Hence
for each $1\leq i\leq k$, there is a finite collection ${\cal Z}_i$ of conjugacy classes of maximal
cyclic subgroups of $G_i$ with the property that  for any $g\in G'-G_i$,
the intersection $gG_jg^{-1}\cap G_i$ is an element of ${\cal Z}_i$.
We now construct a double\footnote{If we double $G'$ without ``twisting''
the edge inclusions then the resulting group $\hat G$ is not
hyperbolic.  But it acts on a  $CAT(0)$ space $X$ so that
$Homeo(\geo X)$ contains $\hat G$ as a finite
index subgroup.}  $G'$ along the collection of subgroups $G_i\defeq Stab(H_i)$, 
$1\leq i\leq k$ as follows: construct a graph of groups ${\cal G}$
with two vertices $v_1,\, v_2$ and $k$ edges $e_1,\ldots,e_k$, where 
$G_{v_i}$ is isomorphic to $G'$ and $G_{e_i}$ is isomorphic to $G_i$.  
Identify $G_{v_i}$ with $G'$.  We choose the embeddings 
$\iota_{ij}:G_{e_i}\ra G_{v_j}$  
so that the image coincides with $G_i\subset G'$, but 
so that the $\iota_{ij}$'s satisfying the following condition:
\begin{quotation}
 (Twisting)\hspace{.5in}  $\iota_{i1}^{-1}({\cal Z}_i)\cap 
\iota_{i2}^{-1}({\cal Z}_i)=\emptyset$.
\end{quotation}

Let $\hat G\defeq\pi_1({\cal G})$,  let
$T$ be the Bass-Serre tree associated with ${\cal G}$,
and let $V$ and $E$ denote the collections of vertices
and edges in $T$ respectively.
$\hat G$ acts (discretely, cocompactly) on a tree of spaces $X$ constructed 
as in section 
\ref{Doubling}, with vertex spaces $X_v,\,v\in V$ and edge
spaces $X_e,\,e\in E$.

\begin{lemma}
$\hat G$ is a hyperbolic group.  All vertex and edge 
groups $G_x,\,x\in V\cup T$  are quasi-convex subgroups
of $\hat G$.
\end{lemma}
\proof
By \cite{BF3}, \cite{Swarup2}, \cite{Mitra} 
it suffices to show that there is an
upper bound on the length of essential annuli (see \cite{BF3}, section 1)
in the graph of groups
${\cal G}$.  Or equivalently, we need to show that there is an
upper bound on the length of any segment in $T$ which is fixed
by a nontrivial element $g\in \hat G$.  We claim that if 
$e_1,\,e_2,\,e_3$ are 3 consecutive edges in the tree $T$,
then $G_{e_1}\cap G_{e_2}\cap G_{e_3}$ is trivial;  for the
twisting condition implies that the intersections
$G_{e_1}\cap G_{e_2}$ and $G_{e_2}\cap G_{e_3}$ are cyclic
subgroups of $G_{e_2}$ with trivial intersection.
\qed

\begin{lemma}
\label{xprops}

1. For every vertex $v\in V$,  $\geo X_v\subset \geo X$ is a 2-sphere.

2. For every edge $e\in E$, $\geo X_e\subset \geo X$ is a circle.

3. If $v_1\neq v_2\in V$ then $\geo X_{v_1}\cap \geo X_{v_2}\approx S^1$  
implies that $v_1$ and $v_2$ are the endpoints of an edge $e\in E$,
and $\geo X_{v_1}\cap \geo X_{v_2}=\geo X_e$.

4. $\cup_{v\in V}\,\geo X_v$ is dense in $\geo X$.

5. Pick $e\in E$, and let $T_1,\, T_2\subset T$ be the two subtrees
that one gets by removing the interior of the edge $e$.  Then
$\geo  X-\geo X_e$ has two connected components, namely the closures
of $(\cup_{v\in T_i}\geo X_v)-\geo X_e$ in $\geo X- \geo X_e$ for $i=1,2$.
\end{lemma}
 The proof of the lemma is similar to arguments from section \ref{Doubling},
so we omit it.

\begin{lemma}
If $\ga\subset\geo X$ is homeomorphic to $S^1$ and $\ga$ 
separates $\geo X$, then $\ga=\geo X_e$ for some $e\in E$.
\end{lemma}

\proof We first claim that $\ga\subset \geo X_v$ for some $v\in V$.
Otherwise by Alexander duality $\geo X_v-\ga$ is connected
for every $v\in V$, and $(\geo X_{v_1}\cup\geo X_{v_2})-\ga$
is connected for any pair of adjacent vertices $v_1,\,v_2\in V$.
By induction this implies that $\cup_{v\in V}\geo X_v-\ga$
is connected.  By part 4 of Lemma \ref{xprops} we conclude that
$\geo X-\ga$ is connected, a contradiction.

Hence we may assume that $\ga\subset \geo X_v$ for some $v\in V$.
Suppose $\ga\neq\geo X_e$ for any $e\in E$ adjacent to $v$.  Then
any point $\xi\in \geo X-\ga$ lies in the same component of 
$\geo X-\ga$ as one of the two components of 
$\geo X_v-\ga$.  By Lemma \ref{circlescut} we can find an edge
$e$ adjacent to $v$ so that $\geo X_e$ intersects both
of the components $U_1,\,U_2$ of $\geo X_v-\ga$.  So we may connect
$U_1$ to $U_2$ within $\geo X_{w}-\ga$ where $w$ is the other endpoint
of $e$.  This contradicts the assumption that $\ga$ separates
$\geo X$.
\qed

Thus, any  homeomorphism $f: \geo X\to \geo X$  preserves the 
collection of circles $\{\geo X_e, e\in E\}$. 

Let ${\cal C}$ denote the collection of unordered triples of circles 
$C_i=\geo X_{e_i}, e_i\in E$, which are {\bf in standard position}, i.e. there 
exists a triple $H_1, H_2, H_3\in {\cal H}$ which are in standard position 
and a homeomorphism $f:  \geo H_1 \cup\geo H_2 \cup \geo H_3 
\to C_1 \cup C_2 \cup C_3$ which carries each circle $\geo H_i$ 
 to one of the circles $C_{j(i)}$. We define the incidence relation for 
elements of ${\cal C}$ the same way as before, let $\Gamma({\cal C})$ 
denote the associated incidence graph. Thus ${\cal C}$ contains the  
subsets ${\cal S}_v$ where ${\cal S}_v$ consists of  triples of circles 
in standard position which are contained in $\geo X_v$. 
Then the incidence graph $\Ga({\cal S}_v)$ is isomorphic to the incidence 
graph of ${\cal S}$, thus it is connected (see part 1 of Lemma 
\ref{circlescut}). 
For each vertex $v\in V$ the union of triples of circles 
$\{C_1, C_2, C_3\}\in {\cal S}_v$ is dense in $\geo X_v$. 

\begin{lemma}
The subgraphs $\Ga({\cal S}_v)$ are the connected components of 
$\Gamma({\cal C})$. 
\end{lemma}
\proof It is enough to show that any  $\{C_1, C_2, C_3\}\in {\cal C}$
is contained in $\geo X_v$ for some $v\in T$, since there
is at most one $\geo X_v$ containing any given pair of circles.

Pick $\{C_1, C_2, C_3\}\in {\cal C}$, with $C_i=\geo X_{e_i}$
for $e_i\in E$.  Note that  $d(e_i,e_j)\leq 1$ for
$1\leq i,j\leq 3$ for otherwise we would have $C_i\cap C_j=\emptyset$.
 Also, observe that if  two of the circles
lie in some $\geo X_v$, then the third one must too (because
$|\geo X_e\cap \geo X_v|\leq 2$ unless $\geo X_e\subset\geo X_v$).
Clearly this forces the edges $e_i$ to share a vertex.
\qed 

Define the incidence graph with the vertex set 
$\{geo X_v,\,v\in T\}$, where the vertices $v,w$ 
are connected by an edge if and only if $\geo X_v \cap \geo X_w\approx S^1$.  
Lemma 18 implies that this graph is isomorphic to the tree $T$. 

\begin{proposition}
\label{spherespreserved}
Any homeomorphism $f: \geo X\to \geo X$ preserves the collection of spheres 
$\{\geo X_v, v\in V\}$.   In particular, $f$ induces an isomorphism of the 
tree $T$.
\end{proposition} 
\proof The homeomorphism $f$ induces an automorphism $f_{\#}$ of the graph 
$\Gamma({\cal C})$, thus it preserves its connected components. Therefore 
for each $v\in V$ there is $w=f_{\#}(v)$ such that 
$f_{\#} \Ga({\cal S}_v)= \Ga({\cal S}_w)$. However 
$$
\cup_{T\in {\cal S}_v} C
$$
is dense in $\geo X_v$. Thus $f$ preserves the collection of spheres 
$\{\geo X_v, v\in V\}$. The paragraph preceeding Proposition implies that 
$f$ induces an automorphism of the tree $T$.  \qed 

\begin{theorem}
\label{rigidexamples}
The homeomorphism group of $\geo X$ contains $\widehat{G}$
as a subgroup of finite index.  Therefore $Homeo(\geo X)$ is a
topologically rigid hyperbolic group.
\end{theorem}
\proof
For every $v\in V$, we identify $\geo X_v$ with
$\geo \H^3$ via a homeomorphism which carries the 
collection $\{\geo X_e\mid e\in E,\,v\subset e\}$  to $\geo{\cal H}$;
this homeomorphism is unique up to a M\"obius transformation
by Proposition \ref{stabismobius}.

Suppose $f\in Homeo(\geo X)$ and $f\restr_{\geo X_v}=id\restr_{\geo X_v}$
for some $v\in V$.  Then $f$ fixes $\geo X_e$ pointwise for
every $e\in E$ containing $v$.  Hence if $v'\in V$ is adjacent to
$v$ then $f(\geo X_{v'})=\geo X_{v'}$.  By Proposition \ref{stabismobius}
 $f\restr_{\geo X_{v'}}$ is a M\"obius transformation.  Either
 $f\restr_{\geo X_{v'}}=id\restr_{\geo X_{v'}}$ or 
 $f\restr_{\geo X_{v'}}$ is a reflection.  But condition 3
 on ${\cal H}$ rules out the latter possibility.  Therefore
 by induction we conclude that $f$  fixes $\geo X_w$ for every $w\in V$,
 and so $f=id$.

Pick $v\in T$, and consider the possibilities for
$f\restr_{\geo X_v}$ where $f\in Homeo(\geo X)$.
There are clearly only finitely many such possibilities
up to post-composition with elements of $\hat G$;
therefore by the preceding paragraph $\hat G$ has finite index in 
$Homeo(\geo X)$.
\qed

\begin{figure}[tbh]
\leavevmode
\centerline{\epsfxsize=3in\epsfbox{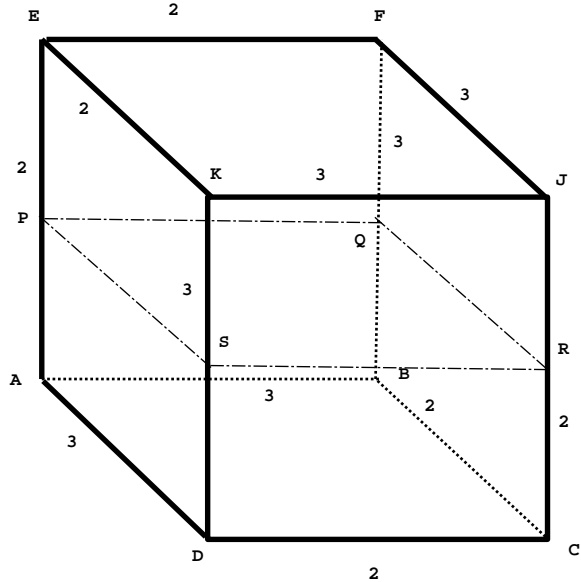}}
\caption{\sl The hyperbolic polyhedron $\Phi$.}
\label{f1}
\end{figure}

\medskip
{\bf An example of a plane configuration ${\cal H}$. }

\begin{figure}[tbh]
\leavevmode
\centerline{\epsfxsize=3in\epsfbox{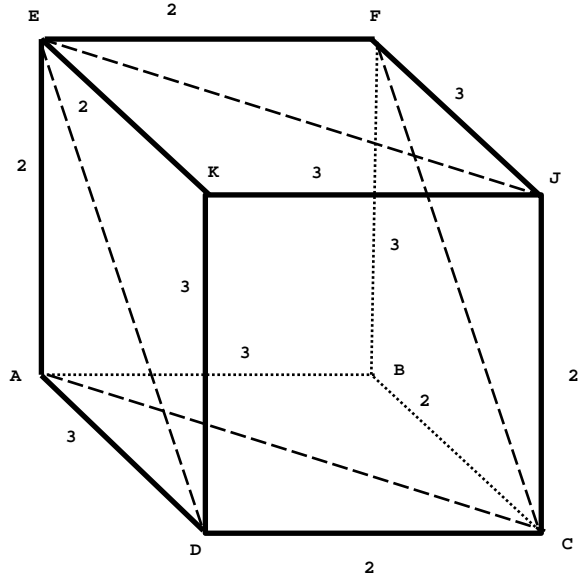}}
\caption{\sl Symmetries of the hyperbolic polyhedron $\Phi$.}
\label{f2}
\end{figure}

\smallskip
We now construct a specific example of a plane configuration ${\cal H}$ 
satisfying the three  
required conditions.  
We start with the 3-dimensional hyperbolic polyhedron $\Phi$ described in  
Figure \ref{f1}: 
the edges of the polyhedron are labelled with $2$ and $3$,  they 
indicate that the corresponding dihedral angles of the polyhedron 
are $\pi/2$ and $\pi/3$ respectively. Such a polyhedron exists by 
Andreev's theorem \cite{Andreev}. Note that $\Phi$ has an order 3 isometry 
$\theta$ 
which is a rotation around the geodesic segment $\overline{CE}$ 
and reflection symmetries in each of  three quadrilaterals, two of which are 
depicted in Figure \ref{f2}. 

The polyhedron $\Phi$ contains three squares which ``bisect'' $\Phi$; 
one of them $\be_1= PQRS$ which is indicated in Figure \ref{f1}, the 
other two $\be_2,\be_3$ are obtained from $\be_1$ by applying the rotation 
$\theta$.

\begin{lem}
The bisectors $\be_1,\be_2,\be_3$ are realized by totally-geodesic 
2-dimensional polygons in $\Phi$ which are orthogonal to the boundary 
of $\Phi$. More precisely, for each $1\le j\le 3$ there is a totally 
geodesic plane $H_j\subset \H^3$  which intersects the same four edges 
of $\Phi$ as $\be_j$ and $H_j$ intersects the faces of $\Phi$ orthogonally.  
\end{lem}  
\proof It is enough to prove the assertion for $\be_1$, the other 
two polygons are obtained via the rotation $\theta$. 
The proof is similar to \cite{Ka}: we first split open the cube 
$\Phi$ combinatorially along the bisector  $\be_1$ into two 
subcubes $\Phi_+$ and $\Phi_-$. Each polyhedron $\Phi_+, \Phi_-$ 
has a face $F_+, F_-$ which corresponds to the bisector  $\be_1$. 
We assign the label $2$ to each edge of $\Phi_{\pm}$ is contained 
in $F_{\pm}$.  Andreev's theorem again implies that $\Phi_+$ 
and $\Phi_-$ can be realized by polyhedra in $\H^3$ (we retain 
the names $\Phi_{\pm}$ for these polyhedra).  Our goal is to 
show that the homeomorphism $F_+\to F_-$ (which is given by 
identification with the bisector  $\be_1$) is isotopic (rel. vertices) 
to an isometry of the hyperbolic polygons. The polyhedron $\Phi$ 
admits a reflection symmetry which fixes the rectangle $EJCA$, 
and this symmetry also acts 
on the polyhedra $\Phi_+, \Phi_-$ and quadrilaterals $F_{\pm}$ so 
that the fixed point sets are the geodesic segments corresponding 
to $\overline{PR}$. However it is clear that there exists a 
unique (up to vertex preserving isotopy) hyperbolic structure 
on quadrilateral $PQRS$ so that the edges are geodesic, angles 
are $\pi/2, \pi/3, \pi/2, \pi/3$ and the quadrilateral has an 
order 2 isometry fixing 
$\overline{PR}$. Thus we have a natural isometry $F_+\to F_-$ 
and we can glue $\Phi_+$  to  $\Phi_-$ using this isometry. The result 
is a hyperbolic polyhedron $\Psi$ which is combinatorially isomorphic 
to $\Phi$ this isomorphism preserves the angles. Thus by uniqueness 
part of Andreev's theorem (alternatively one can use Mostow rigidity 
theorem) the polyhedra $\Phi,\Psi$ are isometric. On the other hand, 
the polyhedron $\Psi$ contains totally geodesic 2-dimensional polygon 
$F_+= F_-$ which is orthogonal to the boundary of $\Psi$. \qed

\begin{figure}[tbh]
\leavevmode
\centerline{\epsfxsize=3in\epsfbox{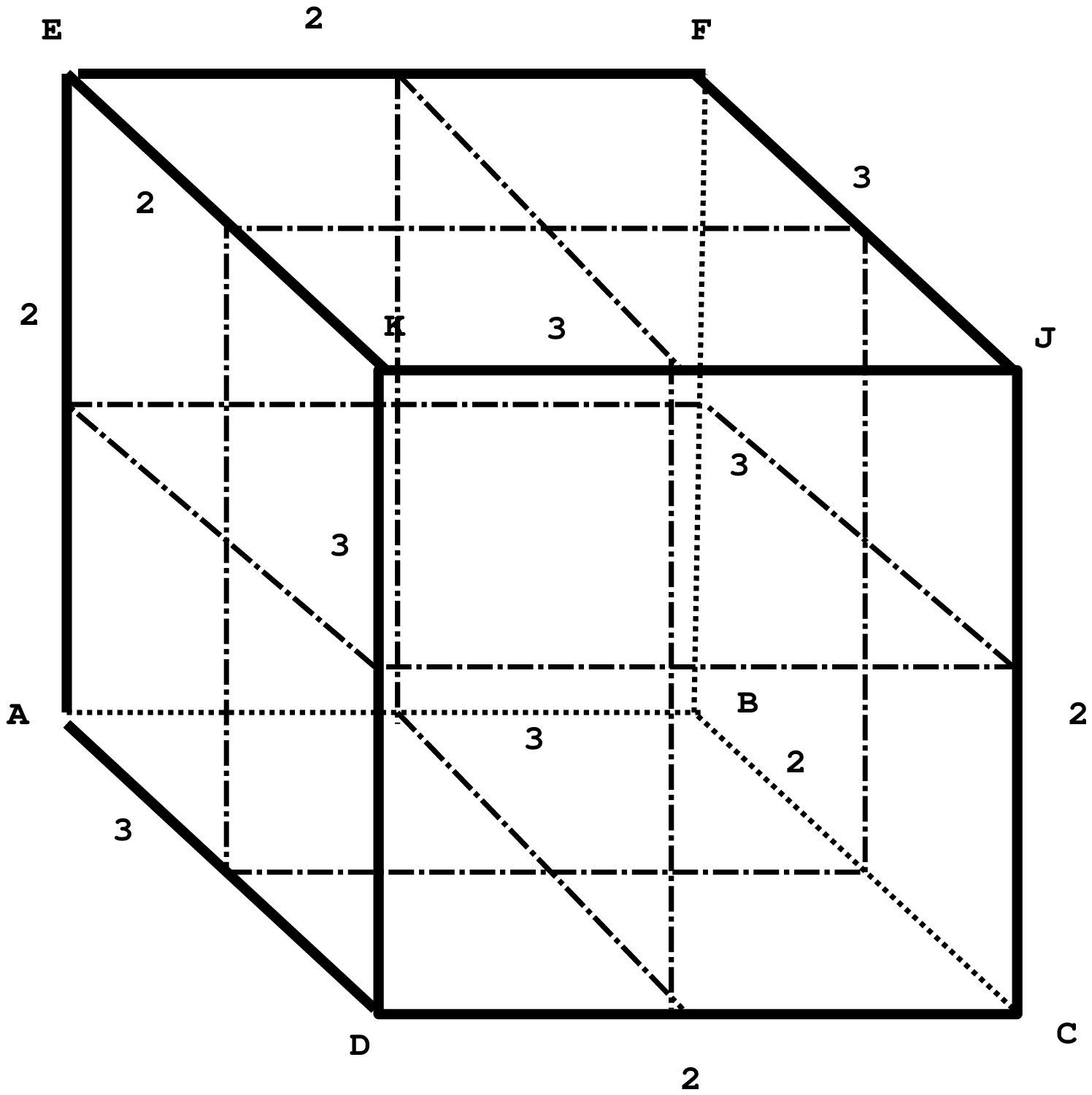}}
\caption{\sl ``Bisectors'' of the hyperbolic polyhedron $\Phi$.}
\label{f3}
\end{figure}

\medskip
We retain the notation $\be_j$ ($j=1,2,3$) for the totally-geodesic 
2-dimensional hyperbolic polygons orthogonal to $\D \Phi$ which realize 
the bisectors $\be_j$. These polygons split $\Phi$ into $8$ subpolyhedra 
$P_i, i=1,...,8$, which are combinatorial cubes. Note that the dihedral 
angles between $\be_j, j=1,2,3$ are all equal and are different from 
$\pi/2$ (otherwise the combinatorial cube  $P_i$ which contains the 
vertex $E$ would have all right angles which is impossible in hyperbolic 
space).   

\begin{figure}[tbh]
\leavevmode
\centerline{\epsfxsize=3in\epsfbox{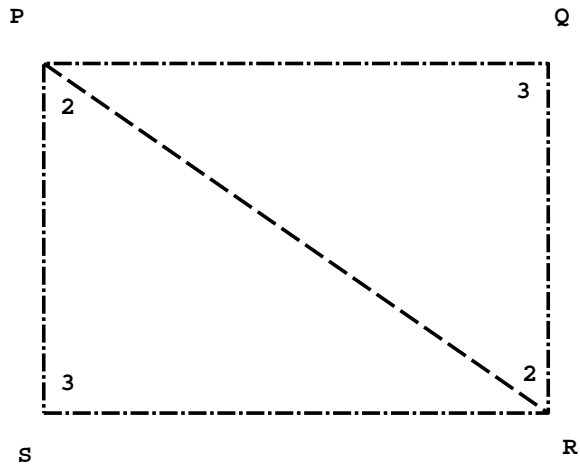}}
\caption{\sl Symmetry of the bisector $\be_1$.}
\end{figure}

Now we construct the collection of planes ${\cal H}$ as follows: let 
${\cal R}\subset Isom(\H^3)$ be the discrete group generated by reflections 
in the faces of $\Phi$; the polyhedron $\Phi$ is a fundamental domain for 
${\cal R}$. The 2-dimensional hyperbolic polygons $\be_j=H_j\cap \Phi$  are 
orthogonal to $\D \Phi$, the plane $H_j$ is invariant under the subgroup 
${\cal R}_j$ of ${\cal R}$ generated by reflections in the faces of $\Phi$ 
which are incident to $\be_j$. The ${\cal R}$-orbit of these hyperplanes 
is ${\cal H}$. Note that 

\smallskip 
(0) If $H$ is a member of ${\cal H}$  and the 
intersection $H\cap \Phi\ne \emptyset$ then $H\cap \Phi$ is equal to 
one of the bisectors $\be_j$. 

\medskip
We next check that ${\cal H}$ satisfies the required properties:

(1) The fundamental domain $\Phi$ for ${\cal R}$ is compact, hence the 
group ${\cal R}$ is a cocompact lattice. 

(2) The complementary regions to ${\cal H}$ in $\H^3$ are finite unions 
of the polyhedra $P_i, i=1,...,8$,  thus they are bounded. 

(3) Let $\rho_j$ be the reflection in the plane $H_j$. Since the planes 
$H_j$, $1\le j\le 3$ are not mutually orthogonal it follows that this 
reflection maps $H_i, i\ne j$, to a plane which does not belong to 
${\cal H}$ (see Property (0) above); it 
follows that $\rho$  does not preserve the configuration ${\cal H}$.

\end{document}